\documentclass[a4paper, twoside]{amsart}
\usepackage[automark]{scrpage2}
\usepackage{hyperref}
\usepackage{amsmath}
\usepackage{amssymb}
\usepackage{amsthm}
\usepackage{amsopn}
\usepackage{amscd}
\usepackage[all]{xypic}
\usepackage{mathtools}
\usepackage{bbm}
\usepackage{listings}
\usepackage{verbatim}
\usepackage{braket}
\usepackage{algorithmic}
\renewcommand{\emph}[1]{{\bf #1}}

\newcommand{\N}{\mathbb{N}}
\newcommand{\Z}{\mathbb{Z}}

\newcommand{\Orbit}{\mathcal{O}}
\newtheorem{definition}{Definition}[section]
\newtheorem{lemma}[definition]{Lemma}
\newtheorem{korollar}[definition]{Corollary}
\newtheorem{satz}[definition]{Theorem}

\DeclareMathOperator{\GL}{GL}

\DeclareMathOperator{\End}{End}
\DeclareMathOperator{\Aut}{Aut}
\DeclareMathOperator{\Hom}{Hom}
\DeclareMathOperator{\Ext}{Ext}

\DeclareMathOperator{\Bild}{Im}

\DeclareMathOperator{\Mod}{\text{mod}}
\DeclareMathOperator{\Mat}{Mat}

\newcommand{\Rep}[1]{\mathrm{Rep}_{#1}}
\newcommand{\dvec}[1]{{\underline{#1}}}
\newcommand{\ledeg}{\le_{\mathrm{deg}}}

\newcommand{\dimv}[1]{\underline{\dim}\left(#1\right)}

\def\clap#1{\hbox to 0pt{\hss#1\hss}}

\numberwithin{equation}{section}
\theoremstyle{remark}

\author{Stefan Wolf}
\date{}
\title{The Hall algebra of a cyclic quiver at $q=0$}
\begin{document}
\bibliographystyle{hplain}
\SelectTips{eu}{10}
\begin{abstract}
We show that the generic Hall algebra of nilpotent representations of an oriented cycle specialised at 
$q=0$ is isomorphic to the generic extension monoid in the sense of Reineke. This continues the work of Reineke in 
\cite{Reineke_genericexts}.
\end{abstract}
\maketitle
\section{The generic extension monoid and the generic Hall algebra}
For the following let $Q$ be a finite quiver with vertex set $Q_0$ and arrow set $Q_1$. We consider only finite 
dimensional and nilpotent representations and modules.For any field $K$ 
we work in $KQ \Mod^0$, the category of nilpotent, finite dimensional $KQ$ representations.
\subsection*{The generic extension monoid}
Fix
an algebraically closed field $K$.
For each dimension vector $\underline{d} \in \N^{Q_0}$ define
\[\Rep{\dvec{d}} := \bigoplus_{\alpha:i \rightarrow j \in Q_1} \Hom_K (K^{d_i}, K^{d_j}).\] 
Each point 
of $\Rep{\underline{d}}$ corresponds to a representation of $Q$ over $K$ and, by choosing bases, for each 
representation $M$ of $Q$ over $K$ there is point $m$ in $\Rep{\dimv{M}}$ such that 
$M \cong m$ as a representation of $Q$.
\[\GL_\dvec{d} := \prod_{i\in Q_0} \GL_{d_i} (K)\]
acts on $\Rep{\dvec{d}}$ by base change. The $\GL_{\dvec{d}}$ orbits are in 1:1 correspondence to
isomorphism classes of $Q$ representations of dimension vector $\dvec{d}$. Denote such an orbit for a finite
dimensional module $M$ by $\Orbit_{M}$. We say that a module $M$ degenerates to $N$, $M \ledeg N$, if
$\Orbit_N \subseteq \overline{\Orbit_M}$, where we take the closure in the Zariski topology.

For two arbitrary sets $U \subseteq \Rep{\dvec{d}}, V \subseteq \Rep{\dvec{e}}$ we define
\begin{align*}
	\mathcal{E}(U,V) := \{ M \in \Rep{\dvec{d} + \dvec{e}} \: |& \: \exists \:A \in U, B \in V \ 
	\text{and a short exact sequence } \\
	 & 0 \rightarrow B \rightarrow M \rightarrow A \rightarrow 0 \}.
\end{align*}

The multiplication on closed irreducible $\GL_\dvec{d}$-stable respectively $\GL_\dvec{e}$-stable 
subvarieties $\mathcal{A} \subseteq \Rep{\dvec{d}}, 
\mathcal{B} \subseteq \Rep{\dvec{e}}$ is defined as:
\[ \mathcal{A} * \mathcal{B} := \mathcal{E}(\mathcal{A}, \mathcal{B})\]
This multiplication is well defined, associative and has a unit: $\Rep{\dvec{0}}$. 
The set of closed irreducible subvarieties of nilpotent representations with 
this multiplication is called the generic extension monoid $\mathcal{M}(Q)$. The composition monoid 
$\mathcal{CM}(Q)$ is the submonoid
generated by the orbits of simple modules without self extensions. For all this see Reineke \cite{reineke_monoid}.
For any word $w = (N_1, \dots, N_r)$ in semisimples we define
$\mathcal{A}_w := \Orbit_{N_1} * \dots * \Orbit_{N_r}$. This is an element of $\mathcal{M}(Q)$.

If $\mathcal{A}$ is a closed irreducible subvariety of some $\Rep{\dvec{d}}$, define
$[\mathcal{A}]:= \mathcal{A}/\GL_{\dvec{d}}$ as a topological space with the quotient topology. 
Hence points of
$[\mathcal{A}]$ correspond to isomorphism classes in $\mathcal{A}$.

\subsection*{The generic Hall algebra}
Now let $K$ be a finite field. Let $M, N, X \in KQ \Mod^0$ be three nilpotent
representations of $Q$ over $K$. Then define
\[ F^X_{M N} := \#\Set{ U \le X | U \text{ submodule}, U \cong N, X/U \cong M } \]
This is a finite number.

Now let $Q$ be Dynkin or an oriented cycle.
Then the isomorphism classes of
nilpotent indecomposable modules are given by some combinatorial set $\Phi$, finitely many for each
dimension vector, independent of the field. Hence an isomorphism class is given by a function 
$\alpha : \Phi \rightarrow \N$ with finite support. For each such $\alpha$ and any field $K$ choose a
module $M(\alpha,K)$
in this isomorphism class. Then there are polynomials $f^{\xi}_{\mu \nu} (q) \in \Z[q]$ such that for each finite
field $K$ we have 
\[ F^{M(\xi,K)}_{M(\mu,K) M(\nu,K)} = f^{\xi}_{\mu \nu} (|K|).\]

Then we define the generic Hall algebra $\mathcal{H}_q (Q)$ to be the free $\Z[q]$ module with basis
$\Set {u_{\alpha} | \alpha \in \N^\Phi \text{ with finite support}}$ and multiplication given by:
\[u_{\mu} \diamond u_{\nu} = \sum_{\xi} f^{\xi}_{\mu \nu} (q) u_{\xi}. \]
The generic composition algebra $\mathcal{C}_q(Q)$ is then the subalgebra of $\mathcal{H}_q(Q)$ generated 
by the simple modules without self extensions. If the quiver is fixed then we often write $\mathcal{H}_q$ and 
$\mathcal{C}_q$ instead of $\mathcal{H}_q(Q)$
and $\mathcal{C}_q(Q)$. Moreover, for convenience we identify for any module $M \cong M(\alpha, K)$, $u_M$ 
with $u_\alpha$ and $f^{X}_{M N}$ with $f^{\xi}_{\mu \nu} (q)$.
For this see Ringel \cite{Ringel_hallalgsandquantumgroups} and \cite{Ringel_hallpolysforrepfiniteheralgs} but
also Hubery \cite{Hubery_hallpolys}.

\section{The composition algebra of a cyclic quiver}
Consider the cyclic quiver $\widetilde{A}_n$ (all arrows in one direction),
$R:=\mathbb{F}_q \widetilde{A}_n$ and $R\Mod^0$,
the category of nilpotent representations over $R$,
where $\mathbb{F}_q$ is a finite field with $q$ elements. 
Let $S_0, S_1, \dots, S_{n-1}$ denote
the simple modules in $R\Mod^0$ numbered in such a way that $\Ext(S_i, S_{i+1})\neq 0$ (now and in the
remainder of this section always count modulo $n$). For generalities on the cyclic quiver see
\cite{dengdu_monbasesforsln}.

Up to isomorphism there is exactly 
one indecomposable representation $S_i[l]$ of length $l$ with socle
$S_i$. Hence the isomorphism classes of modules in $R\Mod^0$ can be described as
elements of $\Pi := \Set{ (\pi^{(0)},\dots,\pi^{(n-1)}) | \pi^{(i)} 
\text{ is a partition } \forall i}$, where each partition $\pi^{(i)}$
describes the indecomposable summands with socle $S_i$. In the notation of the 
previous section $\Phi = \{0, \dots, n-1 \}  \times \N$ and the elements of
$\Pi$ are in obvious bijection to elements of $\N^\Phi$ with finite support. So identify these
functions with $\Pi$.

Now let $M$ be an arbitrary module isomorphic to a $\pi \in \Pi$ and
we denote with $u_M$ or $u_\pi$ its symbol in the generic Hall algebra of $R\Mod^0$. If 
$p=(p_1,p_2,\dots,p_m)$ ($p_1 \ge \dots \ge p_m > 0$) is
a partition we denote by $p(k):=(p_1, \dots, p_{k-1}, p_{k+1},\dots, p_m)$ 
the partition 
obtained by deleting the $k$th component and by $p-1:= (p_1-1,\dots,p_m-1)$. If $T \subseteq \{1,\dots,m\}$ we
denote by $p_T$ the subpartition of $\pi$ by taking the components at the
positions in $T$. For $\pi \in \Pi$ let $\pi_{a,k}$ be equal
to the set of partitions created by deleting the $k$th component of $\pi^{(a-1)}$ and setting
$\pi_{a,k}^{(a)}= \pi^{(a)} \cup (\pi^{(a-1)}_k+1)$ (take the $k$th component of the $(a-1)$st
partition, add one to it and then push it to the right).
\begin{lemma}
	Let $S_i$ be a simple module and $p$ a partition with $n$ components.
	Let $A \subseteq \{1,\dots,m\}$, $|A|=k$ and 
	$B = \{1, \dots, m\} \backslash A$. Let 
	$M:=S_{i-1}[p_A-1] \oplus S_i[p_B]$ and $X=S_i[p]$.
	
	Then $f_{M S_i^k}^{X} = 0 \mod q$ if $p_A \neq p_{\{1,\dots,k\}}$ and
	$f_{M S_i^k}^{X} = 1 \mod q$ if $p_A = p_{\{1,\dots,k\}}$.

	In such a situation define 
	\[Q(X,k) := S_{i-1}[p_{\{1,\dots,k\}}-1] \oplus S_i [p_{\{k+1,\dots,m\}}].\]
	\label{cyclic-one-semisimple-factor}
\end{lemma}
We will see later that $Q(X,k)$ is the quotient of $X$ by $S_i^k$ with the smallest orbit dimension.
\begin{proof}
	Each submodule isomorphic to $S_i^k$ is given by 
	an inclusion of $S_i^k$ into the socle of $X$. Since $S_i^k$ is semisimple and
	$\End(S_i)=\mathbb{F}_q$ these inclusions correspond to elements of
	$\Mat(m,k)$ with rank $k$. Two inclusions $f_1, f_2 \in \Mat(m,k)$
	describe the same submodule if $\Bild(f_1)=\Bild(f_2)$. Column
	transforms on a matrix don't change the image. Using the Gauss
	algorithm we can bring the matrix to column echelon form (this corresponds
	to applying an automorphism on $S_i^k$, $\Aut(S_i^k)=\GL_k(\mathbb{F}_q)$). Hence
	a submodule is given by one $f \in \Mat(m,k)$ of the form:
	\[\left(\begin{array}{@{}c@{}}\xymatrix @=2pt{
	1& 0 \ar@{ *{\cdot} }[rrr]& & & 0 \\
	0& *{*} \ar@{ *{\cdot} }[rrr]& & & 0 \\
	0& 1 \ar@{ *{\cdot} }[rrr]& & & 0 \\
	\\
	\\
	0\ar@{ *{\cdot} }[uuu]& 0 \ar@{ *{\cdot} }[uuu] \ar@{ *{\cdot} }[rrr]& & & 0 \ar@{ *{\cdot} }[uuu] \\
	}\end{array}\right)\]
	Here the $j$th row corresponds to the summand of $X$ isomorphic to $S_i[p_j]$.
	This is actually a certain cell in the cell decomposition of the Grassmannian
	$Gr {m \choose k} $. Now we want to show that the quotient of this submodule
	just depends on the cell.
	
	Let $J$ be the set of indices of the pivot rows, i.e. the rows with a leading $1$. 
	Now fix a $J$ (i.e. fix a cell). Then
	replacing each $*$ with an arbitrary element of $\mathbb{F}_q$ gives us a
	different submodule of $X$.
	
	For all $j < l$ there are embeddings $\phi_{l,j}: S_i[p_l] 
	\rightarrow S_i[p_j]$, unique up to scalar multiply. Without loss of
	generality we can choose them in such a way that
	$\phi_{j,p} \circ \phi_{l,j} = \phi_{l,p}$.	These map the socle of
	one module to the socle of the other one.

	Hence one can do row transforms which add $\lambda$ times a row 
	to a row more to the top by applying an automorphism to $X$, so without
	changing the factor. Hence our matrix becomes one with just ones and
	zeros, at most one $1$ in each row and column. Hence the quotient is
	\[ S_{i-1}[p_J-1] \oplus S_i[p_{\{1,\dots,m\} \backslash J}]\]
	Therefore the quotient is $M$ iff $\pi_J = \pi_A$. Hence
	$J$ determines the quotient (i.e. the quotient is independent of the element of the cell). 
	If $J$ is not completely embedded
	in the top part of the matrix (i.e. $J=\{1,\dots,k\}$), then
	the number of submodules in this cell is $q^r$, where $r$ is the number of $*$ respectively the 
	dimension of the cell. Moreover, if
	$J=\{1,\dots,k\}$ then the number of submodules is $1$. Hence
	the claim follows.
\end{proof}
\begin{korollar}
	Let $X$ be an extension of $M$ by a simple $S_a$, $X$ corresponding
	to $\rho \in \Pi$ and $M$ to $\pi \in \Pi$, where $\rho = \pi_{a,k}$. Let $l= \pi^{(a-1)}_k +1$, 
	$N_j:=\Set{ i | \rho^{(a)}_i = j}$,
	$n_j := \left| N_j \right|$ and 
	$M_j := \Set{ i | \rho^{(a)}_i > j}$,
	$m_j := \left| M_j \right|$.
	% = n_{j+1} - n_j $.
	
	Then \[f_{M S} ^X (q) = \frac{q^{n_l}-1}{q-1} q^{m_l}.\]
\end{korollar}
\begin{proof}
	Using the proof of the previous lemma we just have to count the elements of the cells which
	give us the right quotient. Since we look at one dimensional subspaces there is
	just one column and $J=\{j\}$ consists just of one element. The factor is $M$ iff $j \in N_l$.
	The number of elements in this cell is $q^{j-1}$. So we obtain
	\[
	f_{M S}^X (q)= \sum_{j \in N_l} q^{j-1} = q^{m_l} \sum_{i=0}^{n_l-1} q^i = \frac{q^{n_l} -1}{q-1} q^{m_l}
	\]
\end{proof}
Now we are able to describe the coefficients modulo $q$ for an extension with a semisimple.
\begin{lemma}
	Let $N=\bigoplus_{i=0}^{n-1} S_i^{a_i}$, $a_i \in \N$ be a semisimple module. Let $M\in R\Mod^0$ and 
	$X=\bigoplus_{i=0}^{n-1} S_i[\pi^{(i)}]$, $\pi \in \Pi$ be arbitrary modules. 
	Then
	$f_{M N} ^X = 1 \mod q$ iff 
	\[ M \cong \bigoplus Q(X_i,a_i) =: Q(X,N) \] with $X_i :=S_i[\pi^{(i)}]$ and $f_{M N}^X = 0 \mod q$ otherwise.
	\label{hallnumber-semisimple-arbitrary}
\end{lemma}
\begin{proof}
    Since $\Hom(S_i^{a_i}, S_j[\pi^{(i)}])=0$ for $i \neq j$ every short exact sequence
	$0 \rightarrow N \rightarrow X \rightarrow M \rightarrow 0$ is the direct sum of
	short exact sequences of the form:
	\[
	0 \rightarrow S_i^{a_i} \rightarrow X_i \rightarrow M_i \rightarrow 0
	\]
	for some modules $M_i$ such that $\bigoplus_{i=0}^{n-1} M_i \cong M$. So we have
	%Do I need a sum?\sum_{\substack{(M_0,\dots,M_{n-1}) :\\ \bigoplus M_i \cong M}} 
	\[f_{M N}^X = \sum_{\substack{(M_0,\dots,M_{n-1}) :\\ \bigoplus M_i \cong M}}
	\prod_{i=0}^{n-1} f_{M_i S_i^{a_i}}^{X_i} \]
	where $X_i := S_i[\pi^{(i)}]$. But now $f_{M_i S_i^{a_i}}^{X_i}$ is non-zero modulo $q$ iff
	\[
	M_i \cong S_{i-1}[\pi^{(i)}_{\{1,\dots,a_i\}}-1] \oplus S_{i}[\pi^{(i)}_{\{a_i+1,\dots\}}] = Q(X_i,a_i)\]
	for all $i$ by lemma (\ref{hallnumber-semisimple-arbitrary}). Moreover, the same lemma yields 
	$f_{Q(X_i,a_i) \; S_i^{a_i}}^{X_i} = 1 \mod q$. Hence we are done.
\end{proof}
Now we show that the quotients with coefficients $1 \mod q$ are maximal with respect to the degeneration order.
\begin{lemma}
	Let $S$ be a simple module, $X_1, X_2$ be two arbitrary modules, $f_1 \colon S \rightarrow X_1$
	and $f_2 \colon S \rightarrow X_2$ be two injections. Moreover, let $g \colon X_1 \rightarrow X_2$
	be a morphism such that $g ( f_1 (S)) = f_2 (S)$.

	Then $X_2  \oplus (X_1 /f_1(S)) \ledeg X_1  \oplus (X_2 /f_2(S))$ and this is even an extension 
	degeneration.
	\label{inclusion_hence_degenerates}
\end{lemma}
\begin{proof}
	We construct an extension degeneration. Let 
	$\pi_1 \colon X_1 \rightarrow X_1 /f_1(S)$
	and $\pi_2 \colon X_2 \rightarrow X_2/f_2(S)$ be the canonical projections and let
	$\overline{g} \colon X_1/f_1(S) \rightarrow X_2/f_2(S)$ be the map induced by $\pi_2 g$ (exists since
	$f_1(S) \subseteq \ker (\pi_2 g)$). 
	By construction we have the following commutative diagram:
	\[\begin{CD}
		0 @>>> S @>{f_1}>> X_1 @>{\pi_1}>> X_1/f_1(S) @>>> 0\\
		@. @| @VV{g}V @VV{\bar g}V\\
		0 @>>> S @>{f_2}>> X_2 @>{\pi_2}>>  X_2/f_2(S) @>>> 0.
	\end{CD}\]
	This is a pushout, therefore
	\[ \xymatrix{ 0 \ar[r] & 
	X_1 \ar[r]^(0.3){\left(\begin{smallmatrix}
		g\\
		\pi_1
	\end{smallmatrix}\right)} & 
	X_2 \oplus X_1/f_1(S) \ar[r]^(0.6){\left(\begin{smallmatrix}
		\pi_2 & -\overline{g}
	\end{smallmatrix}\right)}
	& X_2/f_2(S) \ar[r]& 0} \]
	is a short exact sequence, the desired extension degeneration.
%	\begin{itemize}
%		\item 
%			$\left(\begin{smallmatrix}
%				\pi_2 & -\overline{g}
%				\end{smallmatrix}\right)
%			\left(\begin{smallmatrix}
%				g\\
%				\pi_1
%			\end{smallmatrix}\right) = \pi_2 g - \overline{g} \pi_1 = \pi_2 g - \pi_2 g =0$
%		\item Let $m \in X_1$ such that $\left(\begin{smallmatrix}
%				g\\
%				\pi_1
%			\end{smallmatrix}\right)m =0$. Then $m \in f_1(S)$ since $\pi_1 (m) =0$, so there
%			is a $s \in S$ such that $m = f_1(s)$. Now $0 = g(m) = g(f_1(s))$ and therefore 
%			$s=0$ since $gf_1$ is an injection. So $\left(\begin{smallmatrix}
%				g\\
%				\pi_1
%			\end{smallmatrix}\right)$ is injective.
%		\item The second map is surjective since $\pi_2$ is surjective.
%		\item For dimension reasons the image of the first map has to be the
%			kernel of the second map.
%	\end{itemize}
\end{proof}
\begin{korollar}
	Let $M$ be any modules and $h \colon S \rightarrow M$
	be any injection. Let $g \in \End(M)$ be any endomorphism such that $g h \neq 0$.
	Then
	\[ M/h(S) \ledeg M/gh(S). \]
	\label{deg_one_zero}
\end{korollar}
\begin{proof}
From lemma (\ref{inclusion_hence_degenerates}) we have an extension of the form
\[
	0 \rightarrow M \rightarrow M \oplus M/h(S) 
	\rightarrow M/gh(S) \rightarrow 0.
	\] 
By a result of Riedtmann \cite[prop. 4.3]{Riedtmann_degens} this yields that 
\[M/h(S) \ledeg M/gh(S).\]
\end{proof}

\begin{korollar}
	Let $X, M \in R\Mod^0$, $N$ a semisimple module such that there is a short exact sequence
	$0 \rightarrow N \rightarrow X \rightarrow M \rightarrow 0$. Then $M$ degenerates to
	$Q(X,N)$.
	\label{degen-minimal-elts}
\end{korollar}
\begin{proof}
	It is enough to show that this is true for every $X_i$ and $M_i$. So let 
	$X=S_i[p]$, $M=S_{i-1}[p_{A}-1] \oplus S_{i}[p_{B}]$ and $N=S_i^a$. We want to apply the previous 
	corollary to $X$, so we just have to choose an endomorphism $f$ which embeds the module $S_{i}[p_A]$ into
	$S_i[p_{\{1,\dots,k\}}]$. But this is no problem at all since the lengths are bigger to the left. 
	So the corollary follows.
	
%	We prove the
%	lemma by induction on $l:=\left| A \: \Delta \: \{1,\dots,a\} \right|$. If $l=0$ then
%	$M=Q(X,N)$, so we are done. Now let $l>0$. Let $r:=\min\{A \backslash \{1,\dots,a\}\}$ and
%	$s:=\min\{\{1,\dots,a\} \backslash A\}$. Both numbers exist since the symmetric difference
%	is non-zero and both sets are finite and have the same cardinality. By definition $r>s$ and
%	therefore $p_r \le p_s$. With this assumptions we get the following commuting diagram from the
%	Auslander-Reiten-Quiver of $K\widetilde{A_n}$:
%	\[\xymatrix{
%	& S_i[p_s] \ar@{->>}[dr]^{g'}	& \\
%	S_{i}[p_r] \ar@{^{(}->}[ur]^f \ar@{->>}[dr]^{g}& 			& S_{i-1}[p_s-1] \\
%	& S_{i-1} [p_r-1] \ar@{^{(}->}[ur]^{f'}& 
%	}
%	\]
%	Hence we have the following short exact sequence:
%	\[ \xymatrix{
%	0 \ar[r] & S_i[p_r] \ar[r]^-{\left(\begin{smallmatrix} f \\ -g\end{smallmatrix}\right)} & S_i[p_s] \oplus S_{i-1}[p_r-1] 
%		\ar[r]^-{\left(\begin{smallmatrix}
%			g' & f'
%		\end{smallmatrix}\right)}& 
%		 S_i[p_s-1] \ar[r] & 0
%	}
%	\]
%	This means that $S_i[p_s] \oplus S_{i-1}[p_r-1]$ degenerates to $S_i[p_r]\oplus S_{i-1}[p_s-1]$. Set 
%	$A':= A \cup \{s\}\backslash\{r\}$ and $B':= B \cup \{r\} \backslash \{r\}$. Then
%	$S_i[p_{B}] \oplus S_{i-1}[p_A-1]$ degenerates to 
%	$S_i[p_{B'}] \oplus S_{i-1}[p_{A'}-1]$ and we are done by induction since $A' \: \Delta \: \{1,\dots,a\} =
%	(A \: \Delta \: \{1,\dots,a\}) \backslash \{r,s\}$.
\end{proof}
Now we are able to describe the monomial elements of $\mathcal{H}_{\widetilde{A_n}}(q)$, the generic Hall
algebra of the cyclic quiver.
\begin{satz}
	Let $w=(N_1,\dots,N_r)$ be a sequence of semisimple modules. 
%	Let
%	\[\mathcal{I}_w := \Set{ \text{isomorphism classes of } M \in R\Mod^0 | M \text{ has a Loewy series of type w}}.\]
	Then
	\[u_w := u_{N_1} \diamond u_{N_2} \diamond \dots \diamond u_{N_r} = \sum_{M \in [\mathcal{A}_w]} u_{M} \mod q\]
	\label{coeff-compmonoid-areone}
\end{satz}
\begin{proof}
	We proof the theorem by induction on $r$. For $r=1$ the statement is trivial. Now let $w'=wN$ for 
	a semisimple module $N$ and show that
	we have $\sum_{M\in [\mathcal{A}_{w'}]} u_M = u_{w'} = u_w \diamond u_N \mod q$. 
	First, it is obvious that
	every isomorphism class appearing in $u_w \diamond u_N$ with non-zero coefficient is in $[\mathcal{A}_{w'}]$. 
	
	Now let $X$ be in $[\mathcal{A}_{w'}]$.
	We have to show that the coefficient of $u_X$ equals $1$.
	This coefficient is $\sum_{M \in [\mathcal{A}_{w}]} f_{M N}^X$. By lemma (\ref{hallnumber-semisimple-arbitrary})
	this sum is $1 \mod q$ iff $Q(X,N)$ is in $\mathcal{A}_{w}$, and $0$ otherwise. So it remains to show that $Q(X,N) \in \mathcal{A}_{w}$.
	Now, since $X \in \mathcal{A}_{w'}$, we know that there is at least one $M \in \mathcal{A}_{w}$ such that there is a short exact
	sequence $ 0 \rightarrow N \rightarrow X \rightarrow M \rightarrow 0$. Moreover, $\mathcal{A}_w$ is closed
	and via lemma (\ref{degen-minimal-elts})
	we know that $M$ degenerates
	to $Q(X,N)$, hence $Q(X,N) \in \mathcal{A}_{w}$.
\end{proof}
%\begin{satz}
%	Let $w$ be a word of simples. Let 
%	\[I := \Set{ \text{iso-classes of } M \in R\Mod | M \text{ has a composition series of type w}}.\]
%	Then
%	\[u_w := u_{S_{w_1}} * u_{S_{w_2}} * \dots * u_{S_{w_n}} = \sum_{M \in I} u_{M}  \mod q\]
%	where $q= |K|$.
%\end{satz}
%\begin{proof}
%	Proof by induction on the length of w. For $|w| = 1$ this is trivial. Now take a $w$ 
%	with length $n$ and show that $u_w * u_{S_a} = u_{w S_a}$ has the desired form
%	for all $a$. Let $X$ be any module with composition series of type $wS_a$. Hence
%	there is a short exact sequence 
%	\[0 \rightarrow S_a \rightarrow X \rightarrow M \rightarrow 0\]
%	where $M$ is a module with composition series $w$. All possible modules $M$ 
%	are created by $X$ just by removing one $S_a$ in the socle. By the
%	previous lemma there is, up to isomorphism, just one $M$ such that 
%	$q \nmid f_{M S_a}^X$, namely the
%	one where $S_a$ is removed from a summand with maximal length.
%
%	Now by the previous lemma we also know that $f_{M S_a}^X = 1 \mod q$ for this $M$,
%	hence, by induction hypothesis, the coefficient of $X$ in $u_{w S_a}$ is $1$ modulo $q$.
%	Moreover we get all $X$ with this prescribed composition series.
%\end{proof}
By using this we obtain the following result.
\begin{satz}
	The map \[
	\Psi \colon \begin{array}{rcl} \Z\mathcal{M}(\widetilde{A_n}) &\rightarrow& \mathcal{H}_0({\widetilde{A_n}}) \\
		\mathcal{A} &\mapsto& \sum\limits_{M \in [\mathcal{A}]} u_M \end{array}\]
	is an isomorphism  of graded rings.
\end{satz}
\begin{proof}
	Since the isomorphism classes of representations of $\widetilde{A_n}$ are independent of the field,
	the map is well defined. Since $\mathcal{M}(\widetilde{A_n})$ is generated by the sets
	$\Orbit_{N}$ for $N$ semisimple we obtain by theorem (\ref{coeff-compmonoid-areone}) that $\Psi$ is 
	a homomorphism of rings. The images of the basis elements $\Orbit_{N}$ are homogeneous of degree
	$\dimv{N}$ in $\mathcal{H}_0$. Now we show that $\Psi$ is an isomorphism of free $\Z$-modules by
	showing that it is an isomorphism for every degree. 
	The degeneration order is a partial order on the isomorphism classes of modules and we have that
	$\Psi(\overline{\Orbit_{X}}) = u_X + \sum_{X <_{\text{deg}} Y} u_Y$. But all elements of $\mathcal{M}(\widetilde{A_n})$
	are orbit closures and we have that $\Psi$ is 
	injective. Moreover the dimension of the graded components of $\mathcal{M}$ and $\mathcal{H}_0$ are
	both equal to the number of isomorphism classes of this degree. So $\Psi$ is an isomorphism.
\end{proof}
\begin{korollar}
	$\Z \mathcal{CM}(\widetilde{A_n})$ is isomorphic to the generic composition algebra of the 
	cyclic quiver at $q=0$.
\end{korollar}
\begin{proof}
	Everything follows from the theorem since $\Psi$ maps $\Orbit_{S_i}$ to $u_{S_i}$ and these are the generators 
	for $\Z \mathcal{CM}(\widetilde{A_n})$ respectively the composition algebra.
\end{proof}
\bibliography{myrefs}
\end{document}